\documentclass[12pt]{article}
\usepackage{amsmath, amssymb}
\usepackage{graphicx}
\usepackage{color}
\newcommand\real{{\rm I\! R}}
\newcommand\nat{{\rm I\! N}}
\newcommand\cmplx{\;\hbox{\vrule height6.8pt width0.8pt depth-0.1pt
           \kern-3.6pt {\rm C}}}
\newcommand\beweisende{\ \hfill$\Box$\break}

\newtheorem{theorem}{Theorem}[section]
\newtheorem{definition}[theorem]{Definition}

\newtheorem{lemma}[theorem]{Lemma}

\newcommand\proof{{\noindent \it Proof.\ }}
\newcommand\grad{\mathop{\rm grad}\nolimits}
\newcommand\divv{\mathop{\rm div}\nolimits}

\def\bql#1{\begin{equation}\label{#1}}
\def\bq{\begin{equation}}
\def\eq{\end{equation}}
\def\eref#1{{\rm (\ref{#1})}}

\usepackage[bookmarksopen, colorlinks, linkcolor =blue, urlcolor =red, citecolor =red, pagebackref=true, menucolor =blue]{hyperref}

\usepackage{array,multirow}
\usepackage{hhline}
{
\usepackage[table]{xcolor}
\usepackage{authblk}

\begin{document}
\title{On the numerical solution of the elastodynamic problem by a boundary integral equation method}

\author[1]{Roman Chapko\thanks{chapko@lnu.edu.ua}}
\author[2]{Leonidas Mindrinos\thanks{leonidas.mindrinos@univie.ac.at}}
\affil[1]{\small Faculty of Applied Mathematics and Informatics, Ivan Franko National University of Lviv, Ukraine.}
\affil[2]{Computational Science Center, University of Vienna,  Austria.}

 \date{}
 \maketitle

\newcommand\solidrule[1][1cm]{\rule[0.5ex]{#1}{0.15mm}}

 \begin{abstract}
A numerical method for the Dirichlet initial boundary value problem for the elastic equation in the exterior and unbounded region of a smooth closed simply connected 2-dimensional domain, is proposed and investigated. This method is based on a combination of a Laguerre transformation with respect to the time variable and a boundary integral equation approach in the spatial variables. Using the Laguerre transformation in time reduces the time-depended problem to a sequence of stationary boundary value problems, which are solved by a boundary layer approach resulting to a sequence of boundary integral equations of the first kind. The numerical discretization and solution are obtained by a trigonometrical quadrature method. Numerical results are included.   

\vspace*{1cm} {\em Keywords:} elastic equation; initial boundary value problem; Laguerre transformation; fundamental sequence;
single and double layer potentials; boundary integral equations of the first kind; trigonometrical quadrature method.
\end{abstract}

\section{Introduction}

The problem of numerically solving time-dependent boundary value problems has a long history. The simplest approach consists of using the FDM, which has the obvious limitations of simple domains and interior problems. The most used general scheme, for this kind of problems, first reduces the dimensions of the problem by some semi-discretization approach and then solves the simpler problem by a numerical method. For example, using the Galerkin method with respect to the spatial variables, as a semi-discretization technique, results to a  Cauchy problem for  a system of ordinary differential equations. Then, single- or multi-steps methods can be applied. 

On the other hand, we can apply semi-discretization with respect to the time variable and reduce the given problem to a set of stationary problems. This can be done by an integral transformation or by Rothe's method. Then, we can apply a numerical method suitable for stationary problems, for example FEM, FDM, BEM and others.

The main aim and result of this paper is to present an integral equation method for the time-dependent elastic equation in an unbounded two-dimensional domain. The possible variants to use an integral equation approach for time-dependent boundary value problems are discussed in \cite {Co}.  If the given differential equation has a fundamental solution, the problem can be reduced to a time-boundary integral equation by direct or indirect methods (see \cite{An} for the elastic equation). Also by second semi-discretization approach the method of integral equation can be applied for receiving stationary problems. 

Our goal is to extend the idea of the Fourier-Laguerre expansion of the solution and consider it as a semi-discretization approach with respect to time \cite{Ch,Ch1,ChJo,ChKr1,ChKr2, GaCh}. Then, we need to solve a recurrent sequence of stationary problems for Navier equations. The advantage of this semi-discretization approach, compared to Laplace transform, is the simple representation of the numerical solution. In our case, it is a partial sum of the Fourier-Laguerre series.  Since the solution domain is unbounded, the use of classical integral equation methods with a volume potential would be ineffective. Instead, we present the solution in terms of a boundary potentials using explicit fundamental solutions of the obtained sequence. This approach gives us the possibility to reduce the stationary differential problems to a sequence of {\it boundary} integral equations. 

The outline of the paper is as follows: In Section 2, we describe the semi-discretization procedure in time via the Laguerre transformation. Then, the initial boundary value problem for the elastic equation is transformed to a sequence of boundary value problems for Navier equations. In Section 3, we show the reduction of the sequence of stationary problems to integral equations merely involving boundary potentials. To do so, we created fundamental solutions for the sequence of stationary equations. In Section 4, we outline and describe how the well-established numerical methods based on trigonometrical quadratures can be adjusted and applied to solve numerically the derived sequence of boundary integral equations.
In the last Section 5, we demonstrate the feasibility of our approach through three numerical examples.

Before closing this section we formulate the problem to be studied. 
 Let $D\subset\real^2$ be a unbounded domain such that its complement is
 bounded and simply connected and  assume that the boundary $\Gamma$
 of $D$ is of class $C^2$. Consider the initial boundary value problem
 for the hyperbolic elastic equation 
\bql{h_eq}
\frac{{\partial^2 u}}{{\partial \,t^2}} = \Delta^* u
\quad \text{\rm in}\quad
\;D \times \left({0,\infty} \right),
\eq
with the Lam\'e operator defined by $\Delta^* :=c_s^2\Delta+ (c_p^2-c_s^2) \nabla  \nabla \cdot $ , supplied with the homogeneous initial conditions 
\bql{h_ic}
\frac{\partial u}{\partial t}(\cdot,0)=u(\cdot,0) = 0 \quad \text{\rm in}\quad  D,
\eq
and the boundary condition
\bql{h_bc}
u= f,
\quad  \text{\rm on} \quad \Gamma \times (0,\infty),
\eq
where $f$ is a given function, satisfying the compatibility
 condition
 $$
  f(x,0)=\frac{\partial f}{\partial t}\;(x,0)=0,
 \quad x\in \Gamma.
 $$
Here the velocities $c_s$ and $c_p$ have the following form 
$$
c_s=\sqrt{\frac{\mu}{\rho}}, \qquad c_p=\sqrt{\frac{\lambda+2\mu}{\rho}},
$$
where $\rho$ is the density, and $\lambda$ and $\mu$ are the Lam\'e constants.

Since we have an unbounded solution domain, we specify that at infinity
\bql{h_in}
u(x,t) \to 0,
\quad
|x| \to \infty ,
\eq
uniformly with respect to all directions $\frac{x}{{|x|}}$ and all $t \in [0,\infty)$.  Other types of boundary conditions can potentially be handled as well.

\section{Semi-discretization in the time}
\setcounter{equation}{0}
 \setcounter{theorem}{0}
 
For the semi-discretization with respect to the time variable in the problem \eref{h_eq}--\eref{h_in} we use the Laguerre transformation. 

The Laguerre polynomials have the explicit representation
$$
L_n(t) = \sum_{k=0}^n \binom{n}{k} \frac{(-t)^k}{k!},
$$
for $n=0,1,\ldots$. These polynomials are orthogonal with respect to the standard $L^2$ inner product in $(0,\infty)$ appended with the weight $e^{-t}$. The following recurrence relation holds
\bql{rec}
(n+1)L_{n+1}(t)=(2n+1-t)L_n(t)-nL_{n-1}(t).
\eq
For our purpose though, we need to note the relation  
\begin{equation}\label{L_rel}
L_{n+1}'(t)=L_n'(t) - L_n(t), \quad n=0,1,\ldots,
\end{equation}
which is immediate from the explicit representation of $L_n$. 

A function being square integrable with respect to the above weight $e^{-t}$ over the interval $(0,\infty)$ can then be expanded in a scaled Fourier-Laguerre series 
$$
v(t) = \kappa\sum\limits_{n = 0}^\infty {{v_n}{L_n}(\kappa t)} 
$$
with Fourier-Laguerre coefficients 
$$
v_n = \int\limits_0^\infty e^{- \kappa t} L_n(\kappa t)v(t)dt ,
\quad
n = 0,1,\ldots\,\, .
$$
These two relations are interpreted as the inverse and direct Laguerre transformations, respectively. Here, $v$ is the given function and the sequence $\{v_n\}$ is the image; $\kappa > 0$ is a fixed scaling parameter. In the rest of this work, we use the above scaled convention as the Laguerre transform.

Let $\{v_n'\}$ be the sequence obtained when applying the  Laguerre transformation to the derivative of a sufficiently smooth function $v$ with $v(0)=0$, then
\begin{equation}\label{L_rel_d}
v_n'=\kappa\sum_{m=0}^n v_m
\end{equation}
and for the second derivative of $v$ with $v(0)=v'(0)=0$ holds
 \bql{lag9}
 v_n^{\prime\prime}=\kappa^2\sum_{m=0}^n(n-m+1)v_m.
 \eq

Applying the Laguerre transformation to the problem \eref{h_eq}--\eref{h_in} with respect to the time variable, together with the relation (\ref{L_rel_d}), we obtain the following sequence of stationary boundary value problems 
\begin{subequations}
\begin{alignat}{2}
\Delta^* {u_n} - {\kappa^2}{u_n} &= \sum_{m=0}^{n-1}\beta_{n-m}u_m , \quad &&\mbox{in }  D,  \label{s_e}\\
{u_n} &= {f_n}, \quad &&
 \text{on } \Gamma, \label{s_bc}\\
 {u_n}(x) &\to 0, \quad &&
|x| \to \infty , \label{s_in}
\end{alignat}
\end{subequations}
where $n = 0,1, \ldots$,   $\beta_{n}=\kappa^2 (n+1)$ and   $\{u_n\}$ and $\{f_n\}$ are the Fourier-Laguerre sequences of coefficients of the functions $u$ and $f$, respectively. 

By the maximum principle and induction we have the following uniqueness result.
\begin{theorem} 
The sequence of stationary problems \eref{s_e}--\eref{s_in} has at most one solution. 
\end{theorem}

The function of the form 
\bql{u_l}
u(x,t) = \kappa\sum\limits_{n = 0}^\infty {{u_n}(x){L_n}(\kappa t)},
\eq
with $u_n$ solving~\eref{s_e}--\eref{s_in}, 
is clearly a solution to the initial boundary value problem \eref{h_eq}--\eref{h_in}. On the other hand, assuming that the solution to \eref{h_eq}--\eref{h_in} has the right smoothness properties such that it can be expanded in time in terms of the Laguerre polynomials it follows that the coefficients will form a sequence $\{u_n\}$ and satisfy~\eref{s_e}--\eref{s_in}. Thus we state the following theorem.
\begin{theorem} 
A sufficiently smooth function~\eref{u_l} is the solution of the time-dependent problem \eref{h_eq}--\eref{h_in} if and only if its Fourier-Laguerre coefficients  
$u_n$ for $n = 0,1,\ldots$, solve the sequence of stationary problems \eref{s_e}--\eref{s_in}. 
\end{theorem}

\section{A boundary integral equations approach for the stationary problems}
\setcounter{equation}{0}
 \setcounter{theorem}{0}
 
First we determine a sequence of fundamental solutions for equations 
\eref{s_e}.
 \begin{definition} 
\label{def_1}
 The sequence of $2\times 2$ matrices 
 $\{E_n(x,y)\}\; n=0,1,\ldots$  is  called 
 fundamental solutions
 of the equations  \eref{s_e}  if
 \bql{Fund}
 \Delta^* E_n(x,y)- \sum_{m=0}^n \beta_{n-m}E_m(x,y)=\delta(x-y)I.
 \eq
 Here $I$  is the $2\times 2$ identity matrix,
 $\delta$ denotes the Dirac function
 and the differentiation in \eref{Fund} is taken with respect to $x$.
 \end{definition}

Let's consider the polynomials
 $$
 v_n(\gamma, r)=
 \sum_{m=0}^{\left[\frac{n}{2}\right]}a_{n,2m}(\gamma)r^{2m},\quad
 \quad
 w_n(\gamma, r)=
 \sum_{m=0}^{\left[\frac{n-1}{2}\right]}a_{n,2m+1}(\gamma)r^{2m+1}
 $$
 \noindent
 for $n=0,1,\ldots,N-1$ ($w_0=0$),
 where the coefficients $a_{n,m}$ satisfy the recurrence relations
\begin{equation*}
\begin{aligned}
a_{n,0}(\gamma) &=1,\quad n=0,1,\ldots,N-1, \\
a_{n,n}(\gamma) &=-\frac{\gamma}{n}\;a_{n-1,n-1}(\gamma),\quad n=1,2,\ldots,N-1,
\end{aligned}
\end{equation*}
and 
 $$
 a_{n,m}(\gamma)=\frac{1}{2\gamma m}
 \left\{4\left[\frac{m+1}{2}\right]^2a_{n,m+1}(\gamma)
 -\gamma^2\sum_{k=m-1}^{n-1}(n-k+1) a_{k,m-1}(\gamma)\right\},
 $$
 for $m=n-1,\ldots,1 .$
 Here $[r]$ denotes the integer part of  $r\ge 0$.
Next we introduce the sequence of functions
\bql{Phi_n}
 \Phi_n(\gamma, r)=K_0(\gamma r)\,v_n(\gamma ,r)
 +K_1(\gamma r)\,w_n(\gamma, r)
\eq
where $K_0$ and $K_1$ are the modified Hankel functions of order  zero and one, respectively.
Throughout this paper all
functions and constants with a negative index number are set equal
to zero.
\begin{lemma}
The following formulas hold
\bql{int1}
\int_{r/a}^\infty \frac{e^{-\kappa t}L_n(\kappa t)}{\sqrt{t^2-(r/a)^2}}dt=\Phi_n (\tfrac{\kappa}{a}, r )
\eq
and
\bql{int2}
\int_{r/a}^\infty \frac{e^{-\kappa t}(\kappa t)^2L_n(\kappa t)}{\sqrt{t^2-(r/a)^2}}dt=\sum_{k=-2}^2 \chi_{k,n}\Phi_{n+k} (\tfrac{\kappa}{a},r ),
\eq
where
$\chi_{-2,n}=n(n-1)$, $\chi_{-1,n}=-4n^2$, $\chi_{0,n}=2(3n^2+3n+1)$, $\chi_{1,n}=-4(n+1)^2$ and $\chi_{2,n}=(n+1)(n+2)$.
\end{lemma}
\proof
Let us consider the fundamental solution for the wave equation
$$
G(x,y;t)=\frac{\theta(t-|x-y|/a)}{\sqrt{t^2-(|x-y|/a)^2}},
$$
where $\theta$ is the Heaviside function. Clearly it satisfies the equation
$$
\frac{1}{a^2}\frac{\partial^2 G(x,y;t)}{\partial t^2}-\Delta G(x,y;t)=\delta(x-y)\delta(t).
$$ 
If we apply the Laguerre transformation to this equation with respect to time, we receive the following sequence
\bql{seq}
\Delta G_n (x,y)-\frac{\kappa^2}{a^2} \sum_{m=0}^n (n-m+1) G_m (x,y)=\delta(x-y), \;\; n=0,1,\ldots
\eq
for the Laguerre coefficients
$$
G_n(x,y)=\int_0^\infty G(x,y;t)e^{-\kappa t}L_n(\kappa t)dt.
$$
In \cite{ChKr2} it was found with the reducing of \eref{seq} to ordinary differential equations and using its exact solution that $G_n(x,y)=\Phi_n(\frac{\kappa}{a},|x-y|)$.  Thus the formula \eref{int1} is proved. 

The recurrence relation \eref{rec} gives us the following representation
\begin{equation}\label{t2L}
t^2 L_n (t)=\sum_{k=-2}^2 \chi_{k,n}L_{n+k}(t) .
\end{equation}
This relation together with \eqref{int1} results to the formula \eqref{int2}.
\beweisende

Note here, that from \eref{t2L} it follows $\sum_{k=-2}^2 \chi_{k,n}=0$.
 Let us introduce the notation $J(x)=\frac{xx^\top}{|x|^2}$ for $x\in \real^2\backslash \{0\}$.
\begin{theorem}\label{theo}
The sequence of matrices  
\bql{E_n} 
E_n(x,y)=\Phi_{1,n}(|x-y|)I + \Phi_{2,n}(|x-y|)J(x-y) 
\eq 
 are fundamental solutions of \eref{s_e}. 
\\ Here
\begin{align*}
\Phi_{\ell,n}(r) &=\frac{(-\ell)^{\ell-1}}{\kappa^2 r^2}\sum_{k=-2}^2 \chi_{k,n}\left(\Phi_{n+k}(\tfrac{\kappa}{c_s},r) - \Phi_{n+k}(\tfrac{\kappa}{c_p},r)\right)+\frac{(-1)^{\ell-1}}{c_p^2}\Phi_{n}(\tfrac{\kappa}{c_p},r) \\
&\phantom{=}+\frac{\ell-1}{c_s^2}\Phi_{n}(\tfrac{\kappa}{c_s},r),
\end{align*}
for $\ell=1,2$.
\end{theorem}
\proof 
We consider the fundamental solution of the time-dependent elastodynamic equation \eref{h_eq} (see \cite{An})
\begin{equation*}
\begin{aligned}
E(x,y;t) &=\left( \frac{t^2\theta(t-r/c_s)}{r^2\sqrt{t^2-(r/c_s)^2}}-\frac{(t^2-(r/c_p)^2)\theta(t-r/c_p)}{r^2\sqrt{t^2-(r/c_p)^2}}\right)I
\\ 
&\phantom{=}+\left( \frac{(2t^2-(r/c_p)^2)\theta(t-r/c_p)}{r^2\sqrt{t^2-(r/c_p)^2}}-\frac{(2t^2-(r/c_s)^2)\theta(t-r/c_s)}{r^2\sqrt{t^2-(r/c_s)^2}}\right)J(x-y),
\end{aligned}
\end{equation*}
where $r=|x-y|$.
From the definition \ref{def_1} it is clear that
\bql{EL}
E_n(x,y)=\int_0^\infty E(x,y;t)e^{-\kappa t}L_n(\kappa t)dt.
\eq
Thus the statement of the theorem follows from \eref{EL} with the use of formulas \eref{int1} and \eref{int2}.
\beweisende
Note that the fundamental matrix $E_0$ from \eref{E_n} coincides  with the fundamental matrix for the harmonic elastodynamic equation (see \cite{An}).

Noe we can analyze the singularities in the fundamental matrix. The modified Hankel functions have the following series representations
 \begin{equation}\label{macdonald}
  K_0(z) = - \left(\ln \frac{z}2+C\right)\,I_0(z)
 + S_0(z),
\quad
 K_1(z)= \frac{1}{z}+\left(\ln \frac{z}2 +C\right)\,I_1(z)
+S_1(z)
 \end{equation}
with
$$
I_0(z)=\sum^\infty_{n=0} \;
 \frac{1}{(n!)^2}\,\left(\frac{z}{2}\right)^{2n},
 \quad
 I_1(z)=\sum^\infty_{n=0} \;
 \frac{1}{n!(n+1)!}\,\left(\frac{z}{2}\right)^{2n+1},
 $$
and
$$
S_0(z)=\sum^\infty_{n=1}
 \frac{\psi(n)}{ (n!)^2} \,\left(\frac{z}{2}\right)^{2n}, \quad
S_1(z)=-\frac{1}{2}\sum^\infty_{n=0}
 \frac{\psi(n+1)+\psi(n)}{ n!(n+1)!}\,\left(\frac{z}{2}\right)^{2n+1}.
$$
  Here, we set $\psi(0)=0$,
 $$
 \psi(n)=\sum_{m=1}^n\frac{1}{m}\;, \quad n=1,2,\ldots,
 $$
 and let $C = 0.57721\ldots$ denote Euler's constant. Thus we can rewrite the functions $\Phi_n$ as follows
$$
\Phi_n(\gamma, r)= \phi_n(\gamma,r)\ln r+\varphi_n(\gamma,r), \quad n=0,1,\ldots,
$$
where
$$
\phi_n(\gamma,r)=-I_0(\gamma r)v_n(\gamma, r)+I_1(\gamma r)w_n(\gamma, r)
$$
and 
$$
\varphi_n(\gamma, r)=[-(C+\ln\frac{\gamma}{2})I_0(\gamma r)+S_0(\gamma r)]v_n(\gamma, r)+[\frac{1}{\gamma r}+(C+\ln\frac{\gamma}{2})I_1(\gamma r)+S_1(\gamma r)]w_n(\gamma, r).
$$
Clearly we have the following asymptotic behavior with respect to $r$
\bql{asym}
\begin{array}{c}
\displaystyle{
\phi_n(\gamma, r)=\epsilon_{n,0}(\gamma)+\epsilon_{n,2}(\gamma)r^2+O(r^4),}
\\
\displaystyle{
\varphi_n(\gamma, r)=\varepsilon_{n,0}(\gamma)+\varepsilon_{n,2}(\gamma)r^2+O(r^4),
}
\end{array}
\eq
with
$$
\epsilon_{n,0}(\gamma)=-a_{n,0}(\gamma), \quad 
\epsilon_{n,2}(\gamma)=-\frac{\gamma^2}{4}a_{n,0}(\gamma)+\frac{\gamma}{2}a_{n,1}(\gamma)-a_{n,2}(\gamma)
$$
and
$$
\varepsilon_{n,0}(\gamma)=-(C+\ln\frac{\gamma}{2})a_{n,0}(\gamma)+\frac{1}{\gamma}a_{n,1}(\gamma),
$$
$$
\varepsilon_{n,2}(\gamma)=(C+\ln\frac{\gamma}{2})(-\frac{\gamma^2}{4}a_{n,0}(\gamma)+\frac{\gamma}{2}a_{n,1}(\gamma)-a_{n,2}(\gamma))
+\frac{\gamma^2}{4}a_{n,0}-\frac{\gamma}{4}a_{n,1}(\gamma)+\frac{1}{\gamma}a_{n,3}(\gamma).
$$
Then we have the following representation for the functions in \eref{E_n}
\bql{Phi_ln}
\Phi_{\ell,n}(r)=\eta_{\ell,n}(r)\ln r +\xi_{\ell,n}(r),\quad \ell=1,2
\eq
with
\begin{align*}
\eta_{\ell,n}(r) &=\frac{(-\ell)^{\ell-1}}{\kappa^2 r^2}\sum_{k=-2}^2 \chi_{k,n}\left(\phi_{n+k}(\tfrac{\kappa}{c_s},r) - \phi_{n+k}(\tfrac{\kappa}{c_p},r)\right)+
\frac{(-1)^{\ell-1}}{c_p^2}\phi_{n}(\tfrac{\kappa}{c_p},r) \\
&\phantom{=}+\frac{\ell-1}{c_s^2}\phi_{n}(\tfrac{\kappa}{c_s},r) ,
\end{align*}
and
\begin{align*}
\xi_{\ell,n}(r) &=\frac{(-\ell)^{\ell-1}}{\kappa^2 r^2}\sum_{k=-2}^2 \chi_{k,n}\left(\varphi_{n+k}(\tfrac{\kappa}{c_s},r) - \varphi_{n+k}(\tfrac{\kappa}{c_p},r)\right)+
\frac{(-1)^{\ell-1}}{c_p^2}\varphi_{n}(\tfrac{\kappa}{c_p},r)\\
&\phantom{=}
+\frac{\ell-1}{c_s^2}\varphi_{n}(\tfrac{\kappa}{c_s},r).
\end{align*}
Taking into account the definition of the coefficients $\chi_{k,n}$ and $a_{n,m}(\gamma)$ and following \eref{asym} we get the asymptotic expansion
\begin{align*}
\eta_{\ell,n}(r)&=\frac{(-\ell)^{\ell-1}}{\kappa^2}\sum_{k=-2}^2 \chi_{k,n}\left(\epsilon_{n+k,2}(\tfrac{\kappa}{c_s})-\epsilon_{n+k,2}(\tfrac{\kappa}{c_p})\right)+\frac{(-1)^{\ell-1}}{c_p^2}\epsilon_{n,0}(\tfrac{\kappa}{c_p})
\\
&\phantom{=}+\frac{\ell-1}{c_s^2}\epsilon_{n,0}(\tfrac{\kappa}{c_s})+O(r^2)
\end{align*}
and
\begin{align*}
\xi_{\ell,n}(r)&=\frac{(-\ell)^{\ell-1}}{\kappa^2}\sum_{k=-2}^2 \chi_{k,n}\left(\varepsilon_{n+k,2}(\tfrac{\kappa}{c_s})-\varepsilon_{n+k,2}(\tfrac{\kappa}{c_p})\right)+\frac{(-1)^{\ell-1}}{c_p^2}\varepsilon_{n,0}(\tfrac{\kappa}{c_p})
\\
&\phantom{=}+\frac{\ell-1}{c_s^2}\varepsilon_{n,0}(\tfrac{\kappa}{c_s})+O(r^2).
\end{align*}
Thus we are convinced that our fundamental sequence has only the logarithmic singularity.

We shall then construct a solution to the sequence of problems~\eref{s_e}--\eref{s_in}. Let $\{U_n\}$ be a sequence of  single-layer potentials 
\bql{s_p}
{U_n}(x) = \frac{1}{{2\pi}}\sum\limits_{m = 0}^n {\int\limits_\Gamma {{E _{n - m}}(x,y){q_m}(y)\,ds(y)}} ,
\quad
x \in D,
\eq
and $\{V_n\}$  be a sequence of double-layer potentials
\bql{d_p}
{V_n}(x) = \frac{1}{{2\pi}}\sum\limits_{m = 0}^n {\int\limits_\Gamma {{T_yE _{n - m}}(x,y){q_m}(y)\,ds(y)}} ,
\quad
x \in D,
\eq
$n=0,1,\ldots$, where ${q_m} \in C(\Gamma)$ are unknown densities, $\{E_n\}$ is the fundamental sequence \eref{E_n} and $T$ is a tracing operator 
\bql{stress}
Tv =\lambda\divv v\, \nu + 2\mu\,(\nu\cdot \grad)\,v + \mu \divv
(Qv) \,Q\nu 
\eq  
with the unitary matrix 
$$
Q=
\begin{pmatrix}
\phantom{-}0 & 1\\
-1 & 0
\end{pmatrix}.
$$

As follows from the representation of fundamental matrices \eref{E_n}  and the expansion \eref{Phi_ln},  the classical jump and regularity properties of the logarithmic potentials (see \cite{Kr}) can be applied also to the present situation. Hence we have the following transformations into sequences of {\it boundary} integral equations.
\begin{theorem}\label{Sol_Seq_Int_Eq} 
The sequence of single-layer potentials \eref{s_p} is a solution of the sequence of boundary value problems \eref{s_e}--\eref{s_in} provided that their densities satisfy the following sequence of boundary integral equations of the first kind 
\bql{s_bie}
\frac{1}{{2\pi}}\int\limits_\Gamma {{{E _0}(x,y)q_n}(y)\,ds(y) = {f_n}(x) - \frac{1}{{2\pi}}} \sum\limits_{m = 0}^{n - 1} {\int\limits_\Gamma {{{E _{n - m}}(x,y)q_m}(y)\,ds(y)}} ,
\quad
x \in \Gamma ,
\eq
for $n = 0,1,\ldots $ .
\\
The sequence of double-layer potentials \eref{d_p} is a solution of the sequence of boundary value problems \eref{s_e}--\eref{s_in} provided that their densities satisfy the following sequence of boundary integral equations of the second kind 
\begin{align*}
\frac{1}{2}q_n(x)+\frac{1}{2\pi}\int\limits_\Gamma T_yE _0 (x,y)q_n(y)\,ds(y) &= f_n (x) - \frac12\sum\limits_{m = 0}^{n - 1}q_m(x) 
\\
&\phantom{=} -\frac{1}{2\pi} \sum\limits_{m = 0}^{n - 1} \int\limits_\Gamma T_y E _{n - m}(x,y)q_m(y)\,ds(y) ,
\quad
x \in \Gamma ,
\end{align*}
for $n = 0,1,\ldots $ .
\end{theorem}

We proceed by investigating integral equations of the first kind \eref{s_bie}. The case of  the sequence of integral equations of the second kind doesn't contain any principal different.

\begin{theorem}\label{satzigl4}
 For any sequence $f_n$ in $C^{1,\alpha}(\Gamma)$ the system \eref{s_bie}
 possesses a unique solution $q_n$ in $C^{0,\alpha}(\Gamma)$.
 \end{theorem}
\proof 
By standard arguments (see \cite{Kr} for the case of the Laplace equation) it can
be seen that the integral equation with logarithmic singularity 
$$
\frac{1}{2\pi}\int\limits_\Gamma E _0(x,y)q_0(y)\,ds(y) = f_0(x), \quad x\in \Gamma 
$$
has a unique solution $q_0\in C^{0,\alpha}(\Gamma)$ for any $f_0$ in $C^{1,\alpha}(\Gamma)$.  Then the statement
of the theorem follows by induction.
\beweisende

\section{A quadrature  method for full discretization}
\setcounter{equation}{0}
 \setcounter{theorem}{0}
 We assume that the boundary curve $\Gamma $
 is  given
 through
 $$
 \Gamma=\{x(s)=(x_1(s),x_2(s)): 0\leq s \leq 2\pi\},
 $$
 where $x:\real\to \real^2$ is $C^1$ and $2\pi$--periodic with
 $|x^\prime(s)| > 0$ for all $s$, such that the orientation of $\Gamma$
 is counter-clockwise.
 Then we
 transform (\ref{s_bie}) into the parametric form
 \begin{equation}\label{parigl1*}
 \frac{1}{2\pi} \int^{2 \pi}_{0}
 H_0 (s, \tau )\psi_n( \tau )\, d \tau
 =g_n(s)
 -\frac{1}{2\pi} \sum_{m=0}^{n-1}\int^{2 \pi}_{0}
 H_{n-m} (s, \tau )\psi_m( \tau )\, d \tau,
 \quad 0\leq s \leq 2\pi,
 \end{equation}
 where we have set
 $\psi_n(s):= |x^\prime(s)|\,q_n (x(s))$,
 $g_n(s):=f_n(x(s))$
 and where the kernels are given by
 $$
  H_n(s , \tau  ) := E_n(x(s), x(\tau ))
 $$
 for $s \neq \tau$ and $n=0,1,\ldots.$ 

The kernels $H_n$ have  logarithmic singularities
 and can be written in the form
 $$
   H_n(s, \tau ) =
 \ln\left(  \frac{4}e \sin^2 \frac{s- \tau}{2}\right)
 H^1_n (s, \tau )
  + H^2_n (s,\tau ),
 $$
where
$$
H^1_n (s, \tau ):=\frac{1}{2}\left[\eta_{1,n}(|x(s)-x(\tau)|)I+\eta_{2,n}(|x(s)-x(\tau)|)J(x(s)-x(\tau))\right]
$$
and
$$
H^2_n (s, \tau ):= H_n(s, \tau ) -
 \ln\left(  \frac{4}e \sin^2 \frac{s- \tau}{2}\right)
 H^1_n (s, \tau )
$$
with the diagonal terms
$$
H^1_n (s, s )=\frac{1}{2}\left(\eta_{1,n}(0)I+\eta_{2,n}(0)\widetilde J(s,s)\right)
$$
and
$$
H^2_n (s, s )= \frac12 \ln(|x'(s)|^2 e)\left(\eta_{1,n}(0)I+\eta_{2,n}(0)\widetilde J(s,s)\right)+\xi_{1,n}(0)I+\xi_{2,n}(0)\widetilde J(s,s).
$$
Here we used the diagonal  values for the matrix $J,$
$$
\widetilde J(s,s)=
\frac{x^\prime(s) x^\prime(s)^\top}{|x^\prime(s)|^2}.
$$
We choose $M\in\nat$ and an equidistant
 mesh
 by setting
 $
 s_k:=\frac{k\pi}{M}\,,\, k=0,\ldots,2M-1,
 $
 and use the following quadrature rules
 \begin{equation}\label{q1}
 \frac{1}{2\pi}
 \int_0^{2\pi}
 f(\tau)
 \ln\left(  \frac{4}e \sin^2 \frac{s_j- \tau}{2}\right) 
 d\tau
 \approx
 \sum_{k=0}^{2M-1}
 R_{|j-k|}\,f(s_k)
 \end{equation}
 and
 \begin{equation}\label{q3}
 \frac{1}{2\pi}
 \int_0^{2\pi}
 f(\tau)\,
 d\tau
 \approx \frac{1}{2M}
 \sum_{k=0}^{2M-1}
 f(s_k)
 \end{equation}
 with the weights
 $$
 R_j:
 =-\frac{1}{2M} \;
 \left\{ 1-2
 \sum_{m=1}^{M-1} \frac{1}{m} \cos \frac{mj\pi}{M} +\frac{(-1)^j}{M} \right\}, \quad j=0,...,2M-1. 
 $$
 
 These quadratures are obtained by replacing the integrand $f$ by
 its trigonometric interpolation polynomial of degree $M$ with
 respect to the grid points $s_k$, $k=0,\ldots,2M-1.$

 We use the quadrature rules \eref{q1}--\eref{q3} to approximate the
 integrals in the  integral equations \eref{parigl1*} and collocate at the nodal points
 to obtain the sequence of linear systems
 $$
  \sum_{k=0}^{2M-1}
 \left\{
 R_{|j-k|} H_0^1(s_j,s_k)
 +\frac{1}{2M}\;H_0^2(s_j,s_k)
 \right\}\psi_{n,M}(s_k)
 =G_{n,M}(s_j),
 \quad j=0,\ldots,2M-1,
 $$
 which we have to solve for the nodal values $\psi_{n,M}(s_j)$.
 For the right
 hand sides we have
 \bql{rechts}
 G_{n,M}(s_j)=g_n(s_j)
 -
  \sum_{m=0}^{n-1}
  \sum_{k=0}^{2M-1}\left\{
 R_{|j-k|}
 H^1_{n-m} (s_j,s_k) +\frac{1}{2M}\; H^2_{n-m} (s_j,s_k)\right\}
 \psi_{m,M}(s_k ).
 \eq
For a more detailed description of this numerical solution method and an error and convergence analysis based on interpreting the above method as a fully discrete
projection method in a H\"older space setting  and in Sobolev space setting we refer to \cite{Kr}.
In particular, this error analysis implies exponential convergence
$$
\|\psi_{n}-\psi_{n,M}\|_{\infty}\le  C_n e^{-\sigma M}
$$
for some positive constants  $C_n$ and $\sigma,$ provided that the boundary values are also analytic. Of course due to the accumulation of the errors the constants $C_n$ will increase with $n$. 


Given the approximate solution $\psi_{n,M}$ of the integral
 equation \eref{parigl1*}, the approximate solution of the
 initial boundary value problem is obtained by first
 evaluating the parametrized form of the potential \eref{s_p}
 using the trapezoidal rule, that is, by
 \begin{equation}\label{u_tilde}
  \tilde u_{n,M}(x)
 =
 \frac{1}{2M}
 \sum_{m=0}^{n}
  \sum_{k=0}^{2M-1}
 E_{n-m}(x,x(s_k))\psi_{m,N}(s_k),
 \quad x\in D,
 \end{equation}
 and then summing up
 \bql{partreihe}
 u_{N,M}(x,t)=\kappa\sum_{n=0}^{N-1} \tilde u_{n,M}(x)L_n(\kappa t)
 \eq
 according to the series \eref{u_l}.

\begin{figure}[t]
\begin{center}
\includegraphics[scale=0.5]{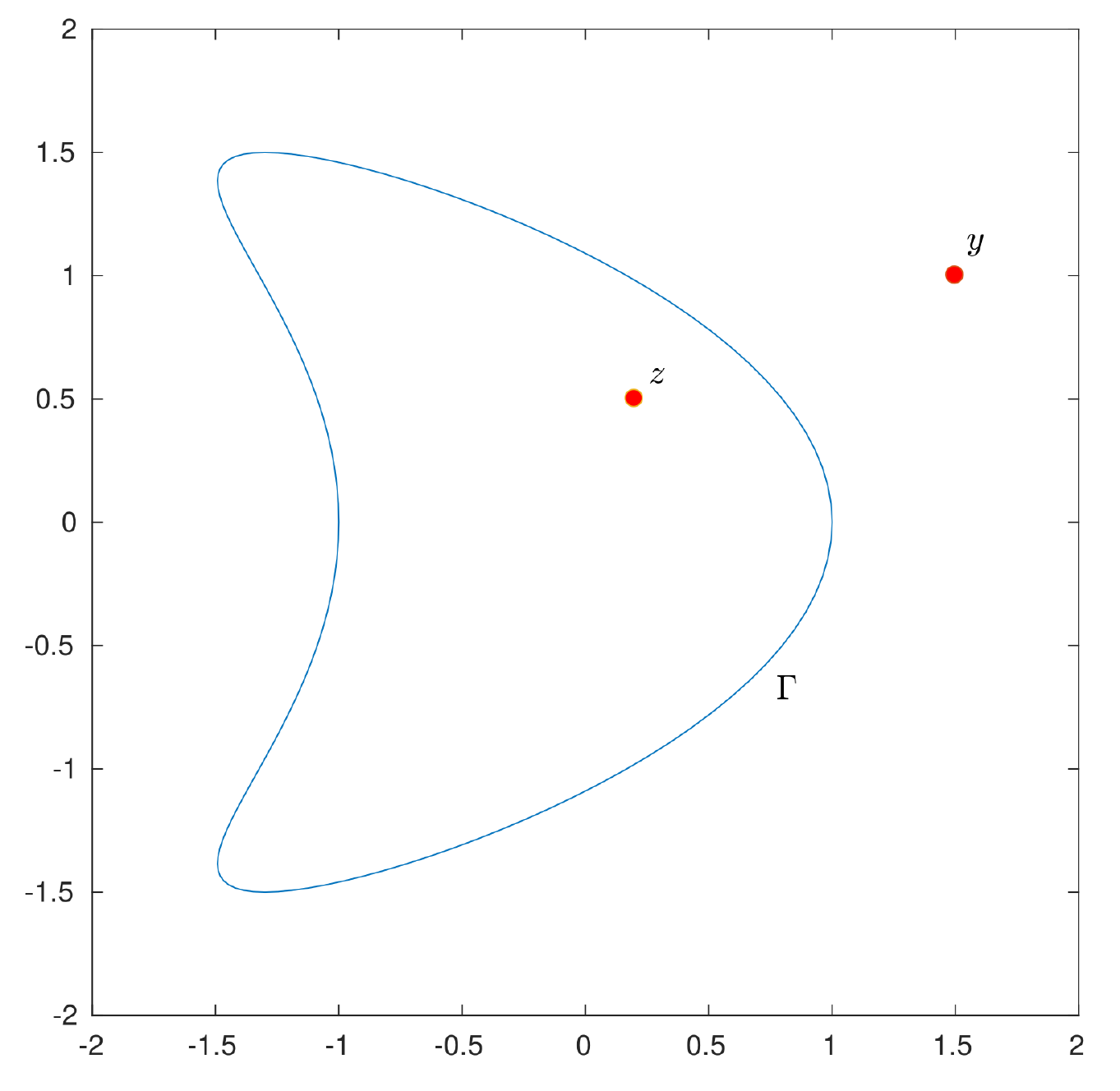}
\caption{The boundary curve $\Gamma$, the source point $z\in \real^2 \setminus D,$ and the measurement point $y\in D.$   }\label{Fig1}
\end{center}
\end{figure}

\section{Numerical results}
\setcounter{equation}{0}
 \setcounter{theorem}{0}

For the numerical examples we consider a kite-shaped boundary with parametrization
\[
x(s) = (\cos(s)+0.65 \cos(2s)-0.65, 1.5 \sin (s) ), \quad s\in [0,2\pi].
\]
In all examples we choose the Lam\'e parameters to be $\lambda = 2, \mu = 1$
and the density $\rho = 1.$ 

\begin{figure}[t]
\begin{center}
\includegraphics[scale=0.7]{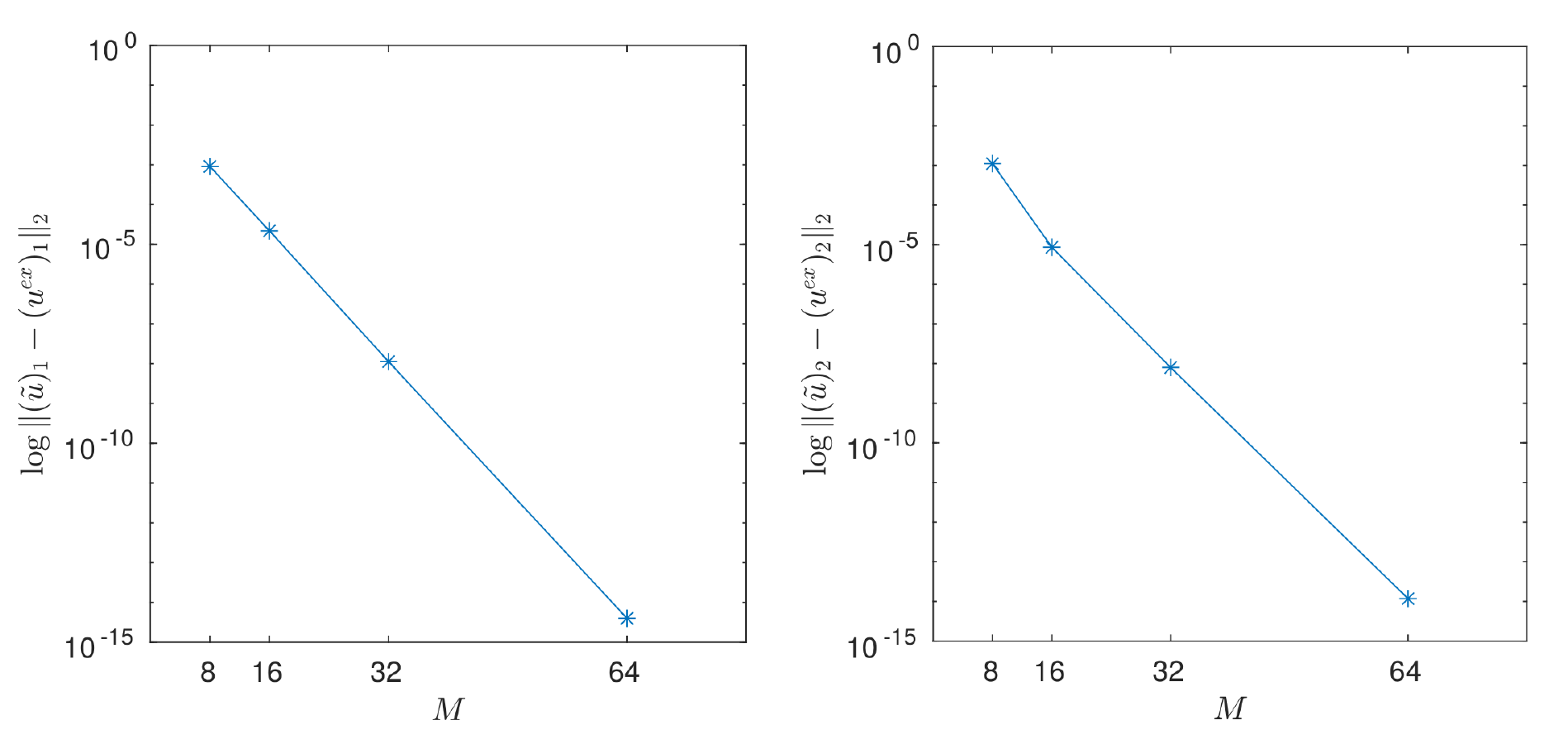}
\caption{The $L^2$ norm of the difference between the computed and the exact solutions in logarithmic scale. In the left picture we see the convergence for the values of Table \ref{table1} and in the right of the Table \ref{table2}.  }\label{Fig2}
\end{center}
\end{figure}

In the first example, we set $\kappa = 1,$ and we test the feasibility of the stationary problems \eref{s_e}--\eref{s_in}. We choose two arbitrary points: a source point $z \in \real^2 \setminus D,$ and a measurement point  $y\in D.$ We define the vector-valued boundary function
\begin{equation}\label{exam_bound1}
f_n(x) = [E_n (x,z)]_1, \quad x\in\Gamma ,
\end{equation}
where $[\cdot]_1,$ denotes the first column of the tensor. Then, the field 
\[
u_n^{ex} (x) := [E_n (x,z)]_1, \quad x\in  D,
\]
is clearly a solution of \eref{s_e}--\eref{s_in} for the boundary function defined above. We consider the points $z = (0.2 , 0.5)$ and $y= (1.5,1),$ see Figure \ref{Fig1}. We present in Tables \ref{table1} and \ref{table2}, the numerical values $\tilde u_{n,M}(y),$ see \eqref{u_tilde}, and compare them with the exact solutions $u_n^{ex} (y),$ for $n=0,1,2$ and varying $M.$ The exponential convergence with respect to the spatial discretization is clearly exhibited, as we can we see also in Figure \ref{Fig2} where we plot the $L^2$ norm of the difference in logarithmic scale.

\begin{table}[t]
\begin{center}
 \begin{tabular}{| c  | c  | c  | c  | } 
 \hline
 $M$ & $(\tilde u_{0,M})_1(y) $ & $(\tilde u_{1,M})_1(y) $ & $(\tilde u_{2,M})_1(y) $  
\\ \hline 
8 & $ 0.293581559232289   $ & $  -0.084483725080856 $ & $  -0.146079079028772 $ \\ 
16 & $ 0.284988364785089   $ & $  -0.092525310787524 $ & $  -0.155666923858005 $ \\
32 & $ 0.285503199323624  $ & $   -0.092138738384510 $ & $  -0.155404881866276 $ \\
64 & $  0.285503741272164  $ & $  -0.092138337605882 $ & $  -0.155404627504594 $ \\
 \cline{1-1}\hhline{~===}

 \multicolumn{1}{c|}{} &  $ (u_0^{ex})_1 (y) $   & $ (u_1^{ex})_1 (y) $ & $ (u_2^{ex})_1 (y) $ \\
 \cline{2-4} 
\multicolumn{1}{c|}{} &  $ 0.285503741272020  $  & $  -0.092138337605708   $  & $-0.155404627504139$ \\ \cline{2-4}
\end{tabular}
\caption{The first components of the computed and the exact solutions of \eref{s_e}--\eref{s_in}, for the specific boundary function \eqref{exam_bound1}, at the measurement point $y=(1.5,1).$}\label{table1}
\end{center}
\end{table}

\begin{table}[!ht]
\begin{center}
 \begin{tabular}{| c  | c  | c  | c  | } 
 \hline
 $M$ & $(\tilde u_{0,M})_2(y) $ & $(\tilde u_{1,M})_2(y) $ & $(\tilde u_{2,M})_2(y) $  
\\ \hline 
8 & 0.081036497084071 &  0.028667287783745 & $-0.011013685642670$ \\ 
16 & 0.071649152048006 &  0.017647071803846 & $-0.021814816353654$ \\
32 & 0.071756738147012  & 0.017837149908881 & $-0.021648258690479$ \\
64 & 0.071756880072043  & 0.017837482038337 & $-0.021647898034988$ \\
 \cline{1-1}\hhline{~===}

 \multicolumn{1}{c|}{} &  $ (u_0^{ex})_2 (y) $   & $ (u_1^{ex})_2 (y) $ & $ (u_2^{ex})_2 (y) $ \\
 \cline{2-4} 
\multicolumn{1}{c|}{} &     0.071756880072350 &  0.017837482039221 & $-0.021647898033812$
 \\ \cline{2-4}
\end{tabular}
\caption{The second components of the computed and the exact solutions of \eref{s_e}--\eref{s_in}, for the specific boundary function \eqref{exam_bound1}, at the measurement point $y=(1.5,1).$}\label{table2}
\end{center}
\end{table}

In the second example we consider the time-dependent problem. We set in \eqref{h_bc} as boundary function
\begin{equation}\label{exam_bound2}
f(x,t) = [E (x,z;t)]_1, \quad x\in \Gamma, \, z\in \real^2 \setminus D, \, t\in (0,\infty).
\end{equation}

Then, the exact solution is given by $u^{ex} (x,t) = [E (x,z;t)]_1, \, x\in  D,$ where its Fourier-Laguerre coefficients satisfy \eref{s_e}--\eref{s_in} for a boundary function $f_n$ as in the first example. We consider different source point $z = (0.4 , 0.2).$ We compare the computed solution $u_{N,M}(y,t),$ formula \eqref{partreihe}, with the exact, considering the truncated form
\begin{equation}\label{truncated}
u^{ex}(N;y,t) = \kappa\sum\limits_{n = 0}^{N-1} {{E_n}(y,z){L_n}(\kappa t)}.
\end{equation}
In Table \ref{table3}, we see the results of the first components, for $\kappa = 1/2,$ at the position $y= (1,1),$ for different values of $M$ and $N$ at various time positions. The values of the second components, at a different position $y=(0.5,-1.5)$ are presented in Table \ref{table4}.

\begin{table}[t]
\begin{center}
 \begin{tabular}{| c | c  | c  | c  | c  | } 
 \hline
$t$ & $M$ & $N=15$ & $N=20 $ & $N=25 $  
\\ \hline 
 \multirow{3}{*}{1} & 32 & 0.556497018896238 & 0.589619531491177 & 0.661666532791282 \\
 & 64 & 0.556495142719243 & 0.589617717103517 & 0.661663974784660 \\ 
 & $(u^{ex})_1$ & \cellcolor[gray]{0.8}0.556495142721411  & \cellcolor[gray]{0.8}0.589617717107183 & \cellcolor[gray]{0.8}0.661663974772766 \\
 \cline{1-1}\hhline{=====}
  \multirow{3}{*}{2} & 32 & 0.512733429371186  & 0.536980500925658 & 0.478568254585875\\
 & 64 &  0.512733615180447 & 0.536980917568962 & 0.478569337106983 \\ 
 & $(u^{ex})_1$ & \cellcolor[gray]{0.8}0.512733615179871  & \cellcolor[gray]{0.8}0.536980917566782 & \cellcolor[gray]{0.8}0.478569337119777  \\
 \cline{1-1}\hhline{=====}
  \multirow{3}{*}{3} & 32 & 0.240000095287178 & 0.141039343209334 & 0.117399715612992 \\
 & 64 & 0.240001133548879 & 0.141039958261864 & 0.117400385116931\\ 
 & $(u^{ex})_1$ & \cellcolor[gray]{0.8}0.240001133547631  & \cellcolor[gray]{0.8}0.141039958259843 & \cellcolor[gray]{0.8}0.117400385110933 \\ \hline
\end{tabular}
\caption{Numerical values of the components of the computed $(u_{N,M})_1$  and the exact solution $(u^{ex}(N))_1$ (rows in grey) of the problem \eref{h_eq}--\eref{h_in} for the boundary function \eqref{exam_bound2}. Here $\kappa = 1/2,$ and $y= (1,1).$}\label{table3}
\end{center}
\end{table}

\begin{table}[!ht]
\begin{center}
 \begin{tabular}{| c | c  | c  | c  | c  | } 
 \hline
$t$ & $M$ & $N=15$ & $N=20 $ & $N=25 $  
\\ \hline 
 \multirow{3}{*}{1} & 32 & $-0.027394214935554$ & $-0.029780601479552$ & $-0.037360920446497$\\
 & 64 & $-0.027394199243065$ & $-0.029780534307372$ & $-0.037360395287884 $ \\ 
 & $(u^{ex})_2$ & \cellcolor[gray]{0.8} $-0.027394199243094$   & \cellcolor[gray]{0.8}$-0.029780534308750$ & \cellcolor[gray]{0.8}$-0.037360395289490$ \\
 \cline{1-1}\hhline{=====}
  \multirow{3}{*}{2} & 32 & $-0.006641666911909$ & $-0.007677796534028$ & $-0.001408341947825$\\
 & 64 &  $-0.006641661641063$  & $-0.007677775150468$ & $-0.001408727582784$ \\ 
 & $(u^{ex})_2$ & \cellcolor[gray]{0.8} $-0.006641661641210$  & \cellcolor[gray]{0.8}$-0.007677775149996$ & \cellcolor[gray]{0.8}$-0.001408727581196$  \\
 \cline{1-1}\hhline{=====}
  \multirow{3}{*}{3} & 32 & 0.008544988585762 & 0.014762278037365 & $0.016882729931766$  \\
 & 64 & 0.008544972022448 & 0.014762136068806 & 0.016882546539029\\ 
 & $(u^{ex})_2$ & \cellcolor[gray]{0.8}0.008544972022701  & \cellcolor[gray]{0.8}0.014762136070815 & \cellcolor[gray]{0.8}$0.016882546538771$ \\ \hline
\end{tabular}
\caption{Numerical values of the components of the computed $(u_{N,M})_2$  and the exact solution $(u^{ex}(N))_2$ (rows in grey) of the problem \eref{h_eq}--\eref{h_in} for the boundary function \eqref{exam_bound2}. Here $\kappa = 1/2,$ and $y= (0.5,-1.5).$}\label{table4}
\end{center}
\end{table}

In the third example, we consider the spatial independent boundary function
\begin{equation}\label{exam_bound3}
f(x,t) = f(t) (1, 1)^\top , \quad \text{for} \quad f(t) =\frac{t^2}4 e^{-t+2},
\end{equation}
which admits the expansion
\[
f(t) = \frac{\kappa e}{4} \sum_{n=0}^\infty \frac{2+\kappa n (\kappa (n-1)-4)}{(\kappa+1)^{n+3}} L_n (\kappa t).
\]
The numerical solution of the problem \eref{h_eq}--\eref{h_in} is presented in Table \ref{table5} (the first component) and in Table \ref{table6} (the second component). Here, we don't know the exact solution but we observe the convergence with respect to the discretization. We set $\kappa =1/2,$ and we compute the solution at the measurement point $y= (0.5,-1.5)$. Again we see the exponential convergence with respect to $M$ and the convergence with respect to $N$.

\begin{table}[t]
\begin{center}
 \begin{tabular}{| c | c  | c  | c  | c  | } 
 \hline
$t$ & $M$ & $N=10$ & $N=15 $ & $N=20 $  
\\ \hline 
 \multirow{3}{*}{1} & 16 &  0.062572635051987 &	0.047833811805535	  &  0.046009612087451
  \\
 & 32 &  0.062555821083946	& 0.047821505207762	 &  0.045995757574454
 \\ 
  & 64 &   0.062555817944294 &	0.047821498423333	&   0.045995752862933
  \\ 
 \cline{1-1}\hhline{=====}
  \multirow{3}{*}{2} & 16 & 0.206527137283687 &	0.230300799019877 &	   0.230584463332204
 \\
 & 32 &  0.206507930918013	& 0.230271954401413	&   0.230552865835900
 \\ 
  & 64 &  0.206507936378995	& 0.230271967614367	&   0.230552881497197
  \\ 
 \cline{1-1}\hhline{=====}
  \multirow{3}{*}{3} & 16 &  0.281765785091472	&   0.290022752671760 &	   0.293418322035422
\\
 & 32 & 0.281753694023185	&   0.290010686252088	&   0.293412971755498
  \\ 
  & 64 & 0.281753701138142	 &  0.290010693385048	 &  0.293412971498783
 \\ \hline
\end{tabular}
\caption{The first component of the numerical solution of the problem \eref{h_eq}--\eref{h_in} for the boundary function \eqref{exam_bound3}, for $\kappa = 1/2,$ at the measurement position $y= (0.5,-1.5).$}\label{table5}
\end{center}
\end{table}

\begin{table}[ht]
\begin{center}
 \begin{tabular}{| c | c  | c  | c  | c  | } 
 \hline
$t$ & $M$ & $N=10$ & $N=15 $ & $N=20 $  
\\ \hline 
 \multirow{3}{*}{1} & 16 & 0.104385591001347	&   0.089274523790871 &	   0.089029847496310
 \\
 & 32 &  0.104414411445436	&   0.089311568889684	 &  0.089061508274427
 \\ 
  & 64 &  0.104414399045771	 &  0.089311556699605	&   0.089061496285464
  \\ 
 \cline{1-1}\hhline{=====}
  \multirow{3}{*}{2} & 16 & 0.242963484605096 &	0.268615650361266 &	   0.270356959200451
 \\
 & 32 & 0.242960849558260 &	0.268601030910481	&   0.270336784091780
  \\ 
  & 64 & 0.242963476811622 	& 0.268601022577740	 &  0.270336775642960
  \\ 
 \cline{1-1}\hhline{=====}
  \multirow{3}{*}{3} & 16 & 0.294596074694631	&   0.301687858830015 &	   0.299978778147128
 \\
 & 32 & 0.294584218004291	 &  0.301675713615639	 &  0.299984785021505
 \\ 
  & 64 & 0.294584215781299	  & 0.301675711488577	&   0.299984782636057
 \\ \hline
\end{tabular}
\caption{The second component of the numerical solution of the problem \eref{h_eq}--\eref{h_in} for the boundary function \eqref{exam_bound3}, for $\kappa = 1/2,$ at position $y= (0.5,-1.5).$}\label{table6}
\end{center}
\end{table}

\end{document}